\documentclass[11pt]{article}
\usepackage{times,mathptm}
\usepackage{amsfonts,amscd,amssymb}
\input diagrams
\newtheorem{definition}{Definition}[section]
\newtheorem{lemma}[definition]{Lemma}
\newtheorem{proposition}[definition]{Proposition}
\newtheorem{corollary}[definition]{Corollary}
\newtheorem{remark}[definition]{Remark}

\newtheorem{theorem}[definition]{Theorem}

\newcommand{\thlabel}[1]{\label{th:#1}}
\newcommand{\thref}[1]{Theorem~\ref{th:#1}}
\newcommand{\selabel}[1]{\label{se:#1}}
\newcommand{\seref}[1]{Section~\ref{se:#1}}
\newcommand{\lelabel}[1]{\label{le:#1}}

\newcommand{\prlabel}[1]{\label{pr:#1}}
\newcommand{\prref}[1]{Proposition~\ref{pr:#1}}
\newcommand{\colabel}[1]{\label{co:#1}}

\newcommand{\relabel}[1]{\label{re:#1}}
\newcommand{\reref}[1]{Remark~\ref{re:#1}}

\newcommand{\eqlabel}[1]{\label{eq:#1}}
\newcommand{\eqref}[1]{(\ref{eq:#1})}

\newenvironment{proof}{{\it Proof.}}{\hfill $ \square $ \vskip 4mm}
\newcommand{\Hom}{{\rm Hom}\,}

\newcommand{\End}{{\rm End}\,}

\def\lan{\langle}
\def\ran{\rangle}

\def\text#1{\mbox{{\rm #1}}}

\def\ol{\overline}
\def\ul{\underline}

\def\leftact{\hbox{$\rightharpoonup$}}
\def\rightact{\hbox{$\leftharpoonup$}}

\def\ot{\otimes}

\def\doublerightleft#1#2{{\lower.2ex\vbox{
\hbox{${\smash{\mathop{\longrightarrow}\limits^{#1}}}$}\vspace*{-4mm}
\hbox{${\smash{\mathop{\longleftarrow}\limits_{#2}}}$}}}}

\begin{document}
\title{Morita Theory for corings and cleft entwining structures\thanks{Research
supported by the bilateral project ``Hopf Algebras in Algebra, Topology,
Geometry and Physics" of the Flemish and
Chinese governments.}}
\author{
S. Caenepeel, J. Vercruysse
\\ Faculty of Applied Sciences\\
Free University of Brussels, VUB\\ B-1050 Brussels, Belgium\and
Shuanhong Wang\thanks{Research partially supported by the Science Foundation
for dinstinguished young scholars
of Henan Province.} \\
Department of Mathematics\\ Henan Normal University\\
Henan, Xinxiang 453002, China}
\date{}
\maketitle

\begin{abstract}
Using the theory of corings, we generalize and unify Morita contexts
introduced by Chase and Sweedler \cite{CS}, Doi \cite{D2}, and Cohen, Fischman
and Montgomery \cite{CFM}. We discuss when the contexts are strict. We apply our
theory corings arising from entwining structures, and this leads
us to the notion of cleft entwining structure.
\end{abstract}

\section{Introduction}\selabel{0}
Let $H$ be a Hopf algebra, $A$ an $H$-comodule algebra, and $B$
the subring of coinvariants. Generalizing a construction due
to Chase and Sweedler \cite{CS}, Doi \cite{D2} gave a Morita
context, connecting $B$ and $\#(H,A)$, and applied this to the
theory of Hopf Galois extensions. In particular, he introduces
the notion of cleft $H$-comodule algebra, and shows that a
cleft $H$-comodule algebra is an $H$-Galois extension.\\
A similar Morita context has been constructed by Cohen, Fischman
and Montgomery in \cite{CFM}. They start from a finite dimensional
Hopf algebra $H$ over a field (or a Frobenius Hopf algebra over a commutative
ring, see \cite{CF}), an $H$-module algebra, and give a Morita context
connecting the smash product $A\# H$ and the ring of invariants.\\
For a finite dimensional Hopf algebra $H$, a left $H$-module algebra
is the same as a right $H^*$-comodule algebra, so it seems obvious
that both contexts then coincide. That this is the case has been pointed
out by Beattie, D\v asc\v alescu and Raianu \cite{BDR}. However, it
is not just a straighforward application of duality principles, since
the connecting bimodules are different in both cases, and since the
Cohen-Fischman-Montgomery structure relies heavily on the fact that
a finite dimensional Hopf algebra is Frobenius (the actions on the
connecting bimodules are defined using the distinguished grouplike).\\
In this paper, we will generalize both contexts. The advantages of our
approach are the following: first, all computations become straightforward
and elementary; secondly, the duality relation between the two contexts and
the connecting bimodules becomes clear, and the r\^ole of Frobenius
type arguments is made clear; in third place, our theory can be applied
in some other particular situations, for examples to generalized smash
products, and to categories of entwined modules; finally, in the infinite
dimensional case, it is clarified why Doi's Morita context is never strict.\\
Our approach is based on a key observation made by Takeuchi \cite{Tak}, that
entwined modules, and, in particular, many kinds of modules such as
relative Hopf modules, Yetter-Drinfeld modules, Doi-Hopf modules etc,
can be viewed as comodules over a certain coring. Takeuchi's observation
has lead to a revived interest in the theory of corings, which goes back
to Sweedler \cite{Sweedler65}. It became clear that corings provide
a unifying and simplifying framework to various topics, such as
Galois theory, descent theory, Frobenius functors and Maschke type Theorems
(see \cite{Br3}, \cite{Br4}, \cite{Gomez}, \cite{Wisbauer}).
Following this philosophy, we can generalize Doi's results, and associate
a Morita context to a coring ${\cal C}$ with a fixed grouplike
element $x$ over a ring $A$ (\seref{2}). In \seref{1a}, we will show
that there is a dual result, which is even more elementary: to a
morphism of rings $i:\ A\to R$, and a right $R$-linear map
$\chi:\ R\to A$ with $\chi(\chi(r)s)=\chi(rs)$ for all $r, s\in R$,
and $\chi(1_R)=1_A$, we can associate a Morita context, which can
in fact be viewed as the Morita context associated to the right $R$-module
$A$, following \cite{Ba}. This Morita context is a generalization
of the Cohen-Fischman-Montgomery context; if $R/A$ is Frobenius, then
both the second connecting bimodule in the context is isomorphic to
$A$ (see \thref{1a.4}). We can give necessary and sufficient conditions
for this Morita context to be strict.\\
To a coring with a fixed grouplike element, we can now associate two
Morita contexts: one to the coring, as mentioned above, and another
one to the dual of the coring, which is a ring. There exists a morphism
between the two contexts, and we have some sufficient conditions for the
two contexts being isomorphic: this is the case when the coring is finitely
generated and projective as an $A$-module, and also when one of the
connecting maps in the Morita context coming from the coring is
surjective, cf. \thref{2.5}.\\
In \seref{3}, we focus attention to the case where the coring ${\cal C}$
arises from an entwining structure $(A,C,\psi)$. We introduce the
notion of cleft entwining structure, and show that cleftness is
equivalent to ${\cal C}$ being Galois in the sense of \cite{Wisbauer},
and $A$ being isomorphic to $A^{{\rm co}{\cal C}}\ot C$ as a left
$A^{{\rm co}{\cal C}}$-module and a right $C$-comodule. The results use
the Morita contexts of the previous Sections. Surprisingly, we were
not able the notion of cleftness to arbitrary corings with a fixed
grouplike element. In Section 5 and 6, we look at factorization 
structures and the smash product, and introduce the notion of
cleft factorization structure.\\
For a coring that is projective, but not necessarily finitely generated,
as an $A$-module, we expect that there is a third Morita context,
connecting the coinvariants and the rational dual of the coring,
generalizing one of the Morita contexts discussed in \cite{BDR}.
This will be discussed in a forthcoming paper.

\section{Preliminaries}\selabel{1}
\subsubsection*{Corings and comodules}
Let $A$ be a ring. The category ${}_A{\cal M}_A$ of $(A,A)$-bimodules
is a monoidal category, and an $A$-coring ${\cal C}$ is a coalgebra in
${}_A{\cal M}_A$, that is
an $(A,A)$-bimodule together with two $(A,A)$-bimodule maps
$$\Delta_{\cal C}:\ {\cal C}\to{\cal C}\ot_A {\cal C} ~~{\rm and}~~
\varepsilon_{\cal C}:\ {\cal C}\to A$$
such that the usual coassociativity and counit properties hold, i.e.
\begin{eqnarray}
(\Delta_{\cal C}\ot_A I_{\cal C})\circ \Delta_{\cal C}&=&
(I_{\cal C}\ot_A\Delta_{\cal C} )\circ \Delta_{\cal C}\eqlabel{1.1.1}\\
(\varepsilon_{\cal C}\ot_A I_{\cal C})\circ \Delta_{\cal C}&=&
(I_{\cal C}\ot_A\varepsilon_{\cal C} )\circ \Delta_{\cal C}
=I_{\cal C}\eqlabel{1.1.2}
\end{eqnarray}
Corings were introduced by Sweedler, see \cite{Sweedler65}.
A right ${\cal C}$-comodule is a right $A$-module $M$ together with
a right $A$-module map $\rho^r:\ M\to M\ot_A{\cal C}$ such that
\begin{eqnarray}
(\rho^r\ot_A I_{\cal C})\circ \rho^r&=&(I_M\ot_A\Delta_{\cal C})\circ \rho^r
\eqlabel{1.1.3}\\
(I_M\ot_A \varepsilon_{\cal C})\circ \rho^r&=& I_M\eqlabel{1.1.4}
\end{eqnarray}
In a similar way, we can define left ${\cal C}$-comodules and
$({\cal C},{\cal C})$-bicomodules. 
We will use the
Sweedler-Heyneman notation for corings and comodules over corings:
$$\Delta_{\cal C}(c)= c_{(1)}\ot_A c_{(2)}~~;~~
\rho^r(m)=m_{[0]}\ot_Am_{[1]}$$
etc. A map $f:\ M\to N$ between (right) ${\cal C}$-comodules is called
a ${\cal C}$-comodule map if $f$ is a right $A$-module map, and
$$\rho^r(f(m))= f(m_{[0]})\ot_Am_{[1]}$$
for all $m\in M$. ${\cal M}^{\cal C}$ is the category of right
${\cal C}$-comodules and ${\cal C}$-comodule maps. In a similar way,
we introduce the categories
$${}^{\cal C}{\cal M},~{}^{\cal C}{\cal M}^{\cal C}, ~{}_{A}{\cal M}^{\cal C}$$
For example, ${}_{A}{\cal M}^{\cal C}$ is the category of right
${\cal C}$-comodules that are also $(A,A)$-bimodules such that the right
${\cal C}$-comodule map is left $A$-linear.\\
Let ${\cal C}$ be an $A$-coring. We write
$${}^*{\cal C}={}_A\Hom({\cal C},A)~~;~~
{\cal C}^*=\Hom_A({\cal C},A)$$
${}^*{\cal C}$ and ${\cal C}^*$ are rings; the multiplication
on ${}^*{\cal C}$ is given by the formula
\begin{equation}\eqlabel{1.1.5}
f\# g=g\circ (I_{\cal C}\ot_A f)\circ \Delta_{\cal C}
\end{equation}
or
\begin{equation}\eqlabel{1.1.6}
(f\# g)(c)=g(c_{(1)}f(c_{(2)}))
\end{equation}
for all left $A$-linear $f,g:\ {\cal C}\to A$ and $c\in {\cal C}$.
The multiplication on ${\cal C}^*$ is given by
$$(f\# g)(c)=f(g(c_{(1)})c_{(2)})$$
The unit is $\varepsilon_{\cal C}$ in both cases. We have
a ring homomorphism
$$i:\ A\to {}^*{\cal C}$$
given by
$$i(a)(c)=\varepsilon_{\cal C}(c)a$$
We easily compute that
\begin{equation}\eqlabel{1.1.7}
(i(a)\#f)(c)=f(ca)~~{\rm and}~~(f\# i(a))(c)=f(c)a
\end{equation}
for all $f\in {}^*{\cal C}$, $a\in A$ and $c\in {\cal C}$.\\
We have a functor
$$F:\ {\cal M}^{\cal C}\to {\cal M}_{{}^*{\cal C}}$$
$F(M)=M$ as a right $A$-module, with right ${}^*{\cal C}$-action
\begin{equation}\eqlabel{1.1.8}
m\cdot f=m_{[0]}f(m_{[1]})
\end{equation}
If ${\cal C}$ is finitely generated and projective as a left
$A$-module, then $F$ is an isomorphism of categories: given
a right ${}^*{\cal C}$-action on $M$, we recover the right
${\cal C}$-coaction by putting
\begin{equation}\eqlabel{1.1.9}
\rho(m)=\sum_j (m\cdot f_j)\ot_A c_j
\end{equation}
where $\{c_j,f_j~|~j=1,\cdots, n\}$ is a finite dual basis of ${\cal C}$
as a left $A$-module.\\
${}^*{\cal C}$ is a right $A$-module, by \eqref{1.1.7}:
$$(f\cdot a)(c)=f(c)a$$
and we can consider the double dual
$$({}^*{\cal C})^*=\Hom_A({}^*{\cal C},A)$$
We have a canonical morphism
$$i:\ {\cal C}\to ({}^*{\cal C})^*;~~i(c)(f)=f(c)$$
We call ${\cal C}$ reflexive (as a left $A$-module) if $i$ is
an isomorphism. If ${\cal C}$ is finitely generated projective as
a left $A$-module, then ${\cal C}$ is reflexive. For any $\varphi\in
({}^*{\cal C})^*$, we then have
$$\varphi=i(\sum_j \varphi(f_j)c_j)$$

\subsubsection*{Galois corings and Descent Theory}
Let ${\cal C}$ be an $A$-coring. Recall that $x\in {\cal C}$
is called grouplike if $\Delta_{\cal C}(x)=x\ot_A x$ and
$\varepsilon_{\cal C}(x)=1$. $G({\cal C})$ is the set of all
grouplike elements in ${\cal C}$. We have the following interpretations
of $G({\cal C})$ (see e.g. \cite[Sec 4.8]{CMZ}, \cite{Br3}).
\begin{eqnarray*}
G({\cal C})
&\cong&\{\rho^r:\ A\to A\ot_A{\cal C}\cong {\cal C}~|~\rho^r~{\rm makes}~A~
{\rm into~a~right}~{\cal C}\hbox{-}{\rm comodule}\}\\
&\cong&\{\rho^l:\ A\to {\cal C}\ot_AA\cong {\cal C}~|~\rho^l~{\rm makes}~A~
{\rm into~a~left}~{\cal C}\hbox{-}{\rm comodule}\}
\end{eqnarray*}
Fix a grouplike element $x$ in ${\cal C}$. We will call $({\cal C},x)$
a coring with fixed grouplike element. 
The associated coactions
on $A$ are given by
$$\rho^r(a)=xa~~;~~\rho^l(a)=ax$$
For a right ${\cal C}$-comodule
$M$, we define the submodule of coinvariants
$$M^{{\rm co}{\cal C}}=\{m\in M~|~\rho(m)=m\ot_A x\}$$
Now let
$$B\subset A^{{\rm co}{\cal C}}={}^{{\rm co}{\cal C}}A=\{b\in A~|~bx=xb\}$$
We have a pair of adjoint functors
$(F,G)$ between the categories ${\cal M}_B$ and ${\cal M}^{\cal C}$,
namely, for $N\in {\cal M}_B$ and $M\in {\cal M}^{\cal C}$,
$$F(N)=N\ot_B A~~{\rm and}~~G(M)=M^{{\rm co}{\cal C}}$$
The unit and counit of the adjunction are
$$\eta_N:\ N\to (N\ot_BA)^{{\rm co}{\cal C}},~~\eta_N(n)=n\ot_B1$$
$$\varepsilon_M:\ M^{{\rm co}{\cal C}}\ot_B A\to M,~~
\varepsilon_M(m\ot_B a)=ma$$
We say that $({\cal C},x)$ satisfies the {\sl Weak Structure Theorem}
if $\varepsilon_M$ is an isomorphism for all $M\in {\cal M}^{\cal C}$,
that is, $G={\bullet}^{{\rm co}{\cal C}}$ is a fully faithful functor.
$({\cal C},x)$ satisfies the {\sl Strong Structure Theorem} if, in addition,
all $\eta_N$ are isomorphisms, or $F$ is fully faithfull, and
therefore $(F,G)$ is an equivalence between categories. Notice that
the Strong Structure Theorem implies that $B=A^{{\rm co}{\cal C}}$.\\
Let $i:\ B\to A$ be a ring homomorphism. It can be verified easily that
${\cal D}=A\ot_BA$, with structure maps
$$\Delta_{\cal D}:\ A\ot_BA\to (A\ot_BA)\ot_A (A\ot_BA)\cong A\ot_BA
\ot_BA~~{\rm and}~~\varepsilon_{\cal D}:\ A\ot_BA\to A$$
given by
$$\Delta_{\cal D}(a\ot_B b)=(a\ot_B 1)\ot_A (1\ot_B b)= a\ot_B 1\ot_B b$$
$$\varepsilon_{\cal D}(a\ot_B b)=ab$$
is an $A$-coring. $1\ot_B 1$ is a grouplike element, and
$({\cal C}, 1\ot_B1)$ is called the {\sl canonical coring}
associated to the ring morphism $i$. Observe that
$${}^*{\cal D}={}_A\Hom(A\ot_B A,A)\cong {}_B\End(A)^{\rm op}$$
If $A$ is finitely generated projective as a left $B$-module, then
${\cal D}$ is reflexive.\\
A right ${\cal D}$-comodule
consists of a right $A$-module $M$ together with a right $A$-module map
$$\rho_M:\ M\to M\ot_A(A\ot_BA)\cong M\ot_B A$$
such that the following
coassociativity and counit condition hold:
\begin{equation}\eqlabel{2.1.1}
m_{[0][0]}\ot_B m_{[0][1]} \ot_B m_{[1]}=m_{[0]}\ot_B 1 \ot_Bm_{[1]}
\end{equation}
and
\begin{equation}\eqlabel{2.1.2}
m_{[0]}m_{[1]}=m
\end{equation}
If $A$ is faithfully flat as a $B$-module, then $({\cal D},1\ot_B1)$
satisfies the Strong Structure Theorem. This was shown in
\cite{Cipolla}; in \cite[Sec. 4.8]{CMZ}, a proof in the coring
language is presented. In fact it is the basic result of descent
theory: an $A$-module $M$ is isomorphic to $N\ot_B A$ for some $B$-module
$N$ if and only if we can define a right ${\cal D}$-coaction on
$A$. In the situation where $A$ and $B$ are commutative, there is an isomorphism
between the category of comodules over the canonical coring, and the
category of {\sl descent data}, as introduced by Knus and Ojanguren in
\cite{KnusO74}, we refer to \cite[Sec. 4.8]{CMZ} for details.
An unpublished result by Journal and Tierney states that, in the
situation where $A$ and $B$ are commutative, $({\cal D},1\ot_B1)$
satisfies the Strong Structure Theorem if and only if
$i:\ A\to B$ is pure as a morphism of $B$-modules. For a proof,
we refer to \cite{Mesablishvili}.\\
Now we return to the general situation, and take an arbitrary
coring $({\cal C},x)$, with fixed grouplike element. Let
$B=A^{{\rm co}{\cal C}}$, and consider the canonical coring
${\cal D}=A\ot_B A$. We have a canonical
coring morphism
$${\rm can}:\ {\cal D}\to {\cal C};~{\rm can}(a\ot_B b)=axb$$
We say that $({\cal C},x)$ is a {\sl Galois coring} if ${\rm can}$
is an isomorphism of corings. In this situation, we obviously have
an isomorphism between the categories ${\cal M}^{{\rm co}{\cal C}}$
and ${\cal M}^{{\rm co}{\cal D}}$.

\begin{proposition}\prlabel{1.1}
Let $({\cal C},x)$ be an $A$-coring with fixed grouplike element,
$B=A^{{\rm co}{\cal C}}$, and ${\cal D}=A\ot_B A$. We then have a
ring homomorphism
$${}^*{\rm can}:\ {}^*{\cal C}\to {}^*{\cal D}\cong {}_B\End(A)^{\rm op};~~
{}^*{\rm can}(f)(a)=f(xa)$$
1) If $({\cal C},x)$ is Galois, then ${}^*{\rm can}$ is an isomorphism.\\
2) If ${}^*{\rm can}$ is an isomorphism, and ${\cal C}$ and ${\cal D}$
are both reflexive (e.g. ${\cal C}$ and $A$ are finitely generated and
projective, resp. as a left $A$-module and a left $B$-module, 
then $({\cal C},x)$ is Galois.\\
3) If $({\cal C},x)$ is Galois, and $({\cal D},x)$ satisfies the
Strong Structure Theorem (e.g. $A$ is faithfully flat as a right
$B$-module), then $({\cal C},x)$ also satisfies the
Strong Structure Theorem.\\
4) If $({\cal C},x)$ satisfies the Weak Structure Theorem, then
$({\cal C},x)$ is Galois.
\end{proposition}

\begin{proof} 1), 2), 3) follow immediately from the observations
made above. ${\cal C}$ is a right ${\cal C}$-module, using $\Delta_{\cal C}$,
and we have a right $B$-module map
$$i:\ A\to {\cal C}^{{\rm co}{\cal C}},~~i(a)=ax$$
It is easily verified that the restriction of $\varepsilon_{\cal C}$
to ${\cal C}^{{\rm co}{\cal C}}$ is an inverse for $i$, so $A$ and
${\cal C}^{{\rm co}{\cal C}}$ are isomorphic in ${\cal M}_B$. Now
$$\varepsilon_{\cal C}={\rm can}:\ A\ot_B A\to {\cal C}$$
is an isomorphism.
\end{proof}

\subsubsection*{Entwined modules}
Let $k$ be a commutative ring, $A$ a $k$-algebra, $C$ a $k$-coalgebra,
and $\psi:\ C\ot A\to A\ot C$ a $k$-linear map satisfying the
following four conditions:
\begin{eqnarray}
&&  (ab)_{\psi}\ot c^{\psi}=  a_{\psi}b_{\Psi}\ot
c^{\psi\Psi}\eqlabel{2.1.1.1}\\
&&(1_A)_{\psi}\ot c^{\psi}=1_A\ot c\eqlabel{2.1.1.2}\\
&&  a_{\psi}\ot \Delta_C(c^{\psi})=
  a_{\psi\Psi}\ot c_{(1)}^{\Psi}\ot c_{(2)}^{\psi}
\eqlabel{2.1.1.3}\\
&&\varepsilon_C(c^{\psi})  a_{\psi} =
\varepsilon_C(c)a\eqlabel{2.1.1.4}
\end{eqnarray}
Here we used the sigma notation
$$\psi(c\ot a)=  a_{\psi}\ot c^{\psi}
=  a_{\Psi}\ot c^{\Psi}$$
We then call $(A,C,\psi)$ a (right-right) entwining structure.
To an entwining structure $(A,C,\psi)$, we can associate an $A$-coring
${\cal C}=A\ot C$. The structure maps are given by the formulas
$$\begin{array}{l}
a'(b\ot c)a=a'ba_{\psi}\ot c^{\psi}\\
\Delta_{\cal C}(a\ot c)=(a\ot c_{(1)})\ot_A (1\ot c_{(2)})\\
\varepsilon_{\cal C}(a\ot c)=a\varepsilon_C(c)
\end{array}$$
An entwined module $M$ is a $k$-module together with a right $A$-action
and a right $C$-coaction, in such a way that
$$\rho^r(ma)=m_{[0]}a_{\psi}\ot m_{[1]}^{\psi}$$
for all $m\in M$ and $a\in A$. The category ${\cal M}(\psi)_A^C$
of entwined modules and $A$-linear $C$-colinear maps is isomorphic to
the category of right ${\cal C}$-comodules.

\subsubsection*{Factorization structures and the smash product}
Let $A$ and $S$ be $k$-algebras, and $R:\ S\ot A\to A\ot S$
a $k$-linear map. We will write
$$R(s\ot a)=a_R\ot s_R=a_r\ot s_r$$
(summation understood). $A\#_R S$ will be the $k$-module
$A\ot S$, with multiplication
\begin{equation}\eqlabel{2.1.2.0}
(a\# s)(b\# t)=ab_R\#s_Rt
\end{equation}
It is straightforward to verify that this multiplication is
associative with unit $1_A\#1_S$ if and only if
\begin{eqnarray}
R(s\ot 1_A)&=& 1_A\ot s\eqlabel{2.1.2.1}\\
R(1_S\ot a)&=& a\ot 1_S\eqlabel{2.1.2.2}\\
R(st\ot a)&=& a_{Rr}\ot s_rt_R\eqlabel{2.1.2.3}\\
R(s\ot ab)&=&a_Rb_r\ot s_{Rr}\eqlabel{2.1.2.4}
\end{eqnarray}
for all $a,b\in A$ and $s,t\in S$. We then call $(A,S,R)$ a
factorization structure, and $A\#_R S$ the smash product
of $A$ and $S$.

\section{The general Morita context}\selabel{1a}
Let $A$ and $R$ be rings, and $i:\ A\to R$ a ring morphism.
We also consider a map $\chi:\ R\to A$ satisfying the following
three conditions, for all $r,s\in R$:
\begin{enumerate}
\item $\chi$ is right $A$-linear;
\item $\chi(\chi(r)s)=\chi(rs)$;
\item $\chi(1_R)=1_A$
\end{enumerate}
It follows from the second condition that $\chi^2=\chi$.
$A$ is a right $R$-module, with structure
$$a\rightact r=\chi(ar)$$
The three conditions on the map $\chi$ can be explained as follows:
$R$ is an algebra in the monoidal category ${}_A{\cal M}_A$
of $A$-bimodules. A map $\chi:\ A\ot_A R=R\to A$ makes $A$
into a right module over this algebra $A$ if and only if it
satisfies this three conditions. This is the dual result of the
fact that grouplike elements on an $A$-coring ${\cal C}$ are
in one-to-one correspondence with right (or left) ${\cal C}$-comodule
structures on $A$.\\
For any right $R$-module $M$, we define
$$M^R=\{m\in M~|~m\cdot r=m\chi(r)\}\cong \Hom_R(A,M)$$
Then $B=A^R=\{b\in A~|~b\chi(r)=\chi(br),~{\rm for~all~}
r\in R\}$ is a subring of $A$,
and
$M^R$ is a right
$B$-module. In fact we obtain a functor
$$G=(\bullet)^R:\ {\cal M}_R\to {\cal M}_B$$
which is a right adjoint of
$$F= \bullet\ot_B A:\ {\cal M}_B\to {\cal M}_R$$
The unit and counit of the adjunction are
$$\eta_N:\ N\to (N\ot_B A)^R,~~\eta_N(n)=n\ot_B 1$$
$$\varepsilon_M:\ M^R\ot_B A\to M,~~M,~~\varepsilon_M(m\ot_B a)=ma$$
In fact, $M^A=\Hom_R(A,M)$. Now consider
\begin{equation}\eqlabel{1a.1.1}
Q=R^R=\{q\in R~|~qr=q\chi(r),~{\rm for~all~}
r\in R\}
\end{equation}
and the map $\varepsilon_Q=\mu:\ Q\ot_B A\to R$,
$\mu(q\ot_B a)=qa$. It is easy to show that $\chi(Q)\subset B$: for all
$q\in Q$ and $r\in R$, we have
$$\chi(q)\chi(r)=\chi(q\chi(r))=\chi(qr)=\chi(\chi(q)r)$$
Recall that $R^*=\Hom_A(R,A)$ is an $(A,R)$-bimodule:\\
$$(a\cdot f\cdot r)(s)=af(rs)$$
for all $a\in A$, $r,s\in S$ and $f\in R^*$.

\begin{lemma}\lelabel{1a.0}
$(R^*)^R\cong A$ as a right $B$-module, and the counit map
$${\rm can}=\varepsilon_{R^*}:\ A\ot_BA\to R^*$$
is given by
$${\rm can}(a\ot_Ba')(r)=a\chi(a'r)$$
for all $a,a'\in A$ and $r\in R$.
\end{lemma}

\begin{proof}
First observe that $f\in (R^*)^R$ if and only if
$$f(rs)=f(\chi(r)s)$$
for all $r,s\in R$. Define
$$j:\ A\to (R^*)^R,~~j(a)(r)=a\chi(r)$$
$$p:\ (R^*)^R\to A:\ p(f)=f(1)$$
It is clear that $j(a)$ is right $R$-linear. Also
$$j(a)(rs)=a\chi(rs)=a\chi(\chi(r)s)=j(a)(\chi(r)s)$$
so $j(a)\in (R^*)^R$. $j$ and $p$ are inverses, since
$$p(i(a))=i(a)(1)=a\chi(1)=a$$
and
$$i(p(f))(r)=f(1)\chi(r)=f(\chi(r))=f(r)$$
Now we compute
$$\varepsilon_{R^*}(a\ot_B a')(r)= (j(a)\cdot a')(r)=
j(a)(a'r)=a\chi(a'r)$$
From this formula, it follows that, for $b\in B$,
$$(j(a)\cdot b)(r)=a\chi(br)=ab\chi(r)=j(ab)(r)$$
so $j$ is right $B$-linear.
\end{proof}

The proof of the following result is now an easy exercise, left
to the reader.

\begin{proposition}\prlabel{1a.1}
With notation as above, $A\in{}_B{\cal M}_R$ and $Q\in{}_R{\cal M}_B$,
and
we have a Morita context $(B,R,A,Q,\tau,\mu)$. The connecting maps
$\mu=\varepsilon_Q:\ Q\ot_B A\to R$ and $\tau:\ A\ot_R Q\to B$ are given by
\begin{equation}
\mu(q\ot_B a)= qa~~{\rm and}~~\tau(a\ot_R q)=a\rightact q=\chi(aq)
\eqlabel{1a.1.2}
\end{equation}
\end{proposition}

\begin{remark}\relabel{1a.1a}
Let $R$ be a ring. Recall from \cite[II.4]{Ba} that we can associate
a Morita context to any right $R$-module $P$. If we consider
$i:\ A\to R$ and $\chi:\ R\to A$ as above, then the Morita context
associated to the right $R$-module $A$ is isomorphic to the Morita
context from \prref{1a.1}. It suffices to observe that
$$B\cong\End_R(A)~~{\rm and}~~Q\cong\Hom_R(A,R)$$
\end{remark}

It is easy to establish when the Morita context is strict.
First let us investigate when $\tau$ is surjective.

\begin{proposition}\prlabel{1a.2}
With notation as in \prref{1a.1}, the following assertions are
equivalent:\\
1) $\tau$ is surjective (and, a fortiori, injective);\\
2) there exists $\Lambda\in Q$ such that $\chi(\Lambda)=1$;\\
3) for all $M\in {\cal M}_R$, the map
$$\omega_M:\ M\ot_R Q\to M^{R},~\omega_M(m\ot_R q)=m\cdot q$$
is an isomorphism;\\
4) $A$ is finitely generated and projective as a right $R$-module.
\end{proposition}

\begin{proof}
$\ul{1)\Rightarrow 2)}$. If $\tau$ is surjective, then there exist
$a_j\in A$ and $q_j\in Q$ such that
$$\tau(\sum_j a_j\ot_Rq_j)=\chi(\sum_j a_jq_j)=1$$
$\Lambda=\sum_j a_jq_j\in Q$, since $Q$ is a left ideal in $R$.\\
$\ul{2)\Rightarrow 3)}$. First observe that $m\cdot q\in M^R$, since
for all $r\in R$:
$$(m\cdot q)\cdot r=m\cdot (qr)=m\cdot q\chi(r)$$
Define
$$\theta_M:\ M^R\to M\ot_R Q,~\eta_M(m)m\ot_R\Lambda$$
We easily compute that
$$\omega_M(\theta_M(m))=m\cdot\Lambda=m\chi(\Lambda)=m$$
for all $m\in M^R$,
and
\begin{eqnarray*}
\theta_M(\omega_M(m\ot_R q))&=& m\cdot q\ot_R\Lambda=m\ot_R q\Lambda\\
&=& m\ot_R q\chi(\Lambda)=m\ot_R q
\end{eqnarray*}
so $\theta_M$ and $\omega_M$ are inverses.\\
$\ul{3)\Rightarrow 1)}$. Observe that $\omega_A=\tau$.\\
$\ul{1)\Rightarrow 4)}$ follows from \cite[Prop. II.4.4]{Ba}, taking
\reref{1a.1a} into account.
\end{proof}

\begin{proposition}\prlabel{1a.2a}
Consider the Morita context from \prref{1a.1}, and assume that
$\tau$ is surjective. Take $\Lambda\in Q$ such that $\chi(\Lambda)=1$.
Then we have the following properties:\\
1) $\Lambda^2=\Lambda$ and $\Lambda R\Lambda=\Lambda B\cong B$.\\
2) The functor $F=\bullet\ot_B A:\ {\cal M}_B\to {\cal M}_R$ is fully
faithful. In other words, for every $N\in {\cal M}_B$, the unit map
$$\eta_N:\ N\to (N\ot_B A)^R,~\eta_N(n)=n\ot 1$$
is fully faithful.\\
3) $B$ is a direct summand of $A$ as a left $B$-module.\\
4) $A$ and $Q$ are generators as left, resp. right $B$-modules.\\
5) $A$ and $Q$ are finitely generated projective as right, resp.
left $R$-modules.\\
6) We have bimodule isomorphisms
$$f:\ A\to {}_R\Hom(Q,R),~~f(a)(q)=qa$$
$$g:\ Q\to \Hom_R(A,R),~~g(q)(a)=qa$$
7) We have algebra isomorphisms
$$k:\ B\to {}_R\End(Q),~~k(b)(q)=qb$$
$$l:\ B\to \End_R(A),~~l(b)(a)=ba$$
\end{proposition}

\begin{proof}
1) Since $\chi(\Lambda)=1$ and $\Lambda\in Q$, we have
$$\Lambda^2=\Lambda\chi(\Lambda)=\Lambda$$
For all $r, s\in R$, we have
$$\chi(r\Lambda)\chi(s)=\chi(r\Lambda\chi(s))=\chi(r\Lambda s)=
\chi(\chi(r\Lambda)s)$$
which implies that $\chi(r\Lambda)\in B$, and
$$\Lambda r\Lambda=\Lambda\chi(r\Lambda)\in\Lambda B$$
so $\Lambda R\Lambda\subset \Lambda B$.\\
Now in the above arguments, take $r=i(b)$, with $b\in B$. It
follows that
$$\Lambda b\Lambda=\Lambda\chi(b\Lambda)=\Lambda b\chi(\Lambda)=
\Lambda b$$
and $\Lambda B\subset \Lambda R\Lambda$. Finally, the right
$B$-module generated by $\Lambda$ is free since $\Lambda b=0$
implies
$$0=\chi(\Lambda b)=\chi(\Lambda)b=b$$
2) If $\tau$ is surjective, then, by standard Morita theory
arguments, the functor $F=\bullet\ot_B A$ is fully faithful,
and has as right adjoint $\bullet\ot_R Q$, and, by the uniqueness
of the adjoint, $\bullet\ot_R Q$ is isomorphic to $G=(\bullet)^R$.
Since $F$ is fully faithful, the unit of the adjunction $(F,G)$
is an isomorphism.\\
3) We define the map ${\rm Tr}:\ A\to B$, ${\rm Tr}(a)=\tau(a\ot_R
\Lambda)=\chi(a\Lambda)$. ${\rm Tr}$ is left $B$-linear, because
$\tau$ is left $B$-linear. ${\rm Tr}$ is a projection, since
$${\rm Tr}(b)=\chi(b\Lambda)=b\chi(\Lambda)=b$$
for all $b\in B$.\\
4)-7) follow from Morita Theory, see \cite[II.3.4]{Ba}.
\end{proof}

We recall from Morita Theory (\cite[II.3.4]{Ba}) that we have
ring morphisms
$$\pi:\ R\to {}_B\End(A)^{\rm op},~~\pi(r)(a)=a\rightact r=\chi(ar)$$
$$\pi':\ R\to \End_B(Q),~~\pi'(r)(q)=rq$$
We also have an $(R,B)$-bimodule map
$$\kappa:\ Q\to {}_B\Hom(A,B),~~\kappa(q)(a)=\chi(aq)$$
and a $(B,R)$-bimodule map
$$\kappa':\ A\to \Hom_B(Q,B),~~\kappa'(a)(q)=\chi(aq)$$
If $\mu$ is surjective, then $\pi$, $\pi'$, $\kappa$ and $\kappa'$
are isomorphisms, and $A$ and $Q$ are finitely generated and projective
as resp. a left and right $B$-module, and a generator as
resp. a right and left $R$-module.

\begin{proposition}\prlabel{1a.3}
Consider the Morita context from \prref{1a.1}. The following
assertions are equivalent:\\
1) $\mu:\ Q\ot_B A\to R$ is surjective;\\
2) the functor $G=(\bullet)^R:\ {\cal M}_R\to {\cal M}_B$
is fully faithful, that is, for all $M\in {\cal M}_R$, the counit
map
$\varepsilon_M:\ M^R\ot_BA\to M$
is an isomorphism;\\
3) $A$ is a right $R$-generator;\\
4) $A$ is projective as a left $B$-module, and $\pi$ is bijective;\\
5) $A$ is projective as a left $B$-module, $\pi$ is injective,
and $\kappa$ is surjective;\\
6) $Q$ is projective as a right $B$-module, $\pi'$ is injective,
and $\kappa'$ is surjective.
\end{proposition}

\begin{proof}
$\ul{1)\Rightarrow 2)}$.
Take $q_j\in Q$ and $a_j\in a$ such that
$\mu(\sum_j q_j\ot_B a_j)=1_R$. For all $m\in M$ and $q\in Q$,
we have that $mq\in M^R$, so we have a well-defined map
$$\theta_M:\ M\to M^R\ot_B A,~~\theta_M(m)=\sum_jmq_j\ot a_j$$
It is clear that $\varepsilon_M\circ \theta_M=I_M$; we also
compute easily that
\begin{eqnarray*}
&&\hspace*{-2cm} \theta_M(\varepsilon_M(m\ot_Ba))=
\theta_M(ma)=\sum_j \sum_jmaq_j\ot a_j\\
&=& \sum_j m\chi(aq_j)\ot_B a_j=\sum_j m\ot_B \chi(aq_j)a_j\\
&=& \sum_j m\ot_B \chi(aq_ja_j)=m\ot_B a
\end{eqnarray*}
$\ul{2)\Rightarrow 1)}$: $\varepsilon_R=\mu$.\\
$\ul{1)\Rightarrow 3,4,5,6)}$: Morita theory (see above).\\
$\ul{1)\Leftrightarrow 3)}$ follows from \cite[Prop. II.4.4]{Ba}, taking
\reref{1a.1a} into account.\\
$\ul{4)\Rightarrow 1)}$: Let $\{a_j,p_j\}$ be a (not necessarily finite)
dual basis of $A$ as a left $B$-module, and put $q_j=\pi^{-1}(p_j)\in R$.
Then
$$\chi(aq_j)=\pi(q_j)(a)=p_j(a)\in B$$
and
$$\pi(q_j\chi(r))(a)=\chi(aq_j\chi(r))=\chi(aq_jr)=\pi(q_jr)(a)$$
$\pi$ is injective, so it follows that
$q_j\chi(r)=q_jr$, and $q_j\in Q=R^R$. Now
$$\mu(\sum_j q_j\ot_Ba_j)=\sum_j q_ja_j=1$$
since $\pi$ is injective and
\begin{eqnarray*}
&&\hspace*{-2cm}\pi(\sum_j q_ja_j)(a)=\sum_j \chi(aq_ja_j)\\
&=&\sum_j \chi(aq_j)a_j=\sum_j p_j(a)a_j=a
\end{eqnarray*}
It follows that $\mu$ is surjective.\\
$\ul{5)\Rightarrow 1)}$: Let $\{a_j,p_j\}$ be a 
dual basis of $A$ as a left $B$-module, and take $q_j\in Q$
such that $\kappa(q_j)=p_j$. Then proceed as in 
$\ul{4)\Rightarrow 1)}$.\\
$\ul{6)\Rightarrow 1)}$: Let $\{q_j,p_j\}$ be a dual basis of
$Q$ as a right $B$-module. We then have, for all $q\in Q$:
Take $a_j\in A$ such that $\kappa'(a_j)=p_j$. Then
$$p_j(q)=\kappa'(a_j)(q)=\chi(a_jq)$$
hence, for all $q\in Q$,
$$q=\sum_j q_jp_j(q)=\sum_j q_j\chi(a_jq)=\sum_j q_ja_jq$$
so
$$\pi'(\sum_j q_jp_j)=\pi'(1_R)$$
and, since $\pi'$ is injective,
$$\mu(\sum_j q_j\ot_Ba_j)=\sum_j q_ja_j=1$$
Therefore $\mu$ is surjective.
\end{proof}

Recall that the ring extension $R/A$ is called {\sl Frobenius}
if there exists an $A$-bimodule map $\ol{\nu}:\ R\to A$ and
$e=e^1\ot_A e^2\in R\ot_AR$ (summation implicitly understood) such that
\begin{equation}\eqlabel{1a.4.1}
re^1\ot_A e^2=e^1\ot_A e^2r
\end{equation}
for all $r\in R$, and
\begin{equation}\eqlabel{1a.4.2}
\ol{\nu}(e^1)e^2=e^1\ol{\nu}(e^2)=1
\end{equation}
This is equivalent to the restrictions of scalars ${\cal M}_R\to
{\cal M}_A$ being Frobenius, which means that its left and right
adjoints are isomorphic (see \cite[Sec. 3.1 and 3.2]{CMZ}).
$(e,\ol{\nu})$ is then called a Frobenius system.

\begin{theorem}\thlabel{1a.4}
Let $i:\ A\to R$ be a morphism of rings, and $\chi:\ R\to A$
a map satisfying the conditions stated at the beginning of this
Section. If $R/A$ is Frobenius, with Frobenius system
$(e,\ol{\nu})$, then $A$ is an $(R,B)$-bimodule, with left $R$-action
$$r\cdot a=\ol{\nu}(ra\chi(e^1)e^2)$$
Then $A\cong Q$ as $(R,B)$-bimodules, and we have a Morita
context
$$(B,R,A,A,\tau,\mu)$$
with connecting maps
$$\mu:\ A\ot_B A\to R:\ \mu(a\ot_B a')=a\chi(e^1)e^2a'$$
$$\tau:\ A\ot_R A\to B:\ \tau(a\ot_R a')=\chi(aa'\chi(e^1)e^2)$$
\end{theorem}

\begin{proof}
Define $\alpha:\ A\to Q$ by $\alpha(a)=a\chi(e^1)e^2$.
$\alpha(a)\in Q$ since
\begin{eqnarray*}
&&\hspace*{-2cm}
\alpha(a)r=a\chi(e^1)e^2r=a\chi(re^1)e^2\\
&=& a\chi(\chi(r)e^1)e^2=a\chi(e^1)e^2\chi(r)=\alpha(a)\chi(r)
\end{eqnarray*}
for all $r\in R$.
The restriction of $\ol{\nu}$ to $Q$ is the inverse of $\alpha$:
\begin{eqnarray*}
&&\hspace*{-2cm}\ol{\nu}(\alpha(a))=
\ol{\nu}(a\chi(e^1)e^2)=a\chi(e^1)\ol{\nu}(e^2)\\
&=&a\chi(e^1\ol{\nu}(e^2))=a\chi(1)=a
\end{eqnarray*}
and
\begin{eqnarray*}
&&\hspace*{-2cm}
\alpha(\ol{\nu}(q))=\ol{\nu}(q)\chi(e^1)e^2=
\ol{\nu}(q\chi(e^1))e^2\\
&=&\ol{\nu}(qe^1)e^2=\ol{\nu}(e^1)e^2q=q
\end{eqnarray*}
for all $a\in A$ and $q\in R$.
$\alpha$ is right $B$-linear, since
$$\alpha(ab)=ab\chi(e^1)e^2=a\chi(be^1)e^2=a\chi(e^1)e^2b=\alpha(a)b$$
for all $a\in A$ and $b\in B$. It is easy to see that the left $R$-action
on $Q$ is transported into the required left $R$-action on $A$.
The rest follows easily from \prref{1a.1}.
\end{proof}

\begin{remark}
Another possible approach to \thref{1a.4} is the following:
if $R/A$ is Frobenius, then $R^*=\Hom_A(R,A)$ and $R$ are isomorphic
as $(A,R)$-bimodules (see \cite[Theorem 28]{CMZ}). Consequently
$$Q=R^R\cong (R^*)^R=A$$
\end{remark}

\section{A Morita context associated to a coring}\selabel{2}
In this Section, $A$ is a ring, ${\cal C}$ is an $A$-coring,
$x\in{\cal C}$ is a fixed grouplike element. Let $R={}^*{\cal C}$
and consider
$$\chi:\ R\to A,~~\chi(f)=f(x)$$
Using \eqref{1.1.7}, we can easily compute that $\chi$ is right
$A$-linear, $\chi(i(\chi(f))\# g)=\chi(f\# g)$, and $\chi(\varepsilon_{\cal C})
=1$. Any right ${\cal C}$-comodule $M$ is also a right
${}^*{\cal C}$-module (see \eqref{1.1.8}), and it is easy to
prove that
$$M^{{\rm co}{\cal C}}\subset M^{{}^*{\cal C}}$$
If ${\cal C}$ is finitely generated and projective as a left $A$-module,
then the converse implication also holds, and the coinvariants coincide
with the invariants. We put
$$B'=A^{{\rm co}{\cal C}}=\{b\in A~|~bx=xb\}
\subset B=A^{{}^*{\cal C}}=\{a\in A~|~f(c)a=af(x),~{\rm for~all~}
c\in {\cal C},~f\in {}^*{\cal C}\}
$$
${}_A\End({\cal C})$ is a left ${}^*{\cal C}$-module: for all
$f\in {}^*{\cal C}$ and $\varphi\in {}_A\End({\cal C})$, 
and $c\in {\cal C}$, we define
$$(f\# \varphi)(c)=\varphi(c_{(1)}f(c_{(2)}))$$
Now let
\begin{eqnarray}
Q'&=&\{q\in {}^*{\cal C}~|~q\#I_{\cal C}=\rho^l\circ q\}\nonumber\\
&=&\{q\in {}^*{\cal C}~|~c_{(1)}q(c_{(2)})=q(c)x,~{\rm for~all~}
c\in {\cal C}\}\eqlabel{2.2.1}
\end{eqnarray}
Observe that
$$Q'\subset Q=({}^*{\cal C})^{{}^*{\cal C}}$$
and $Q'=Q$ if ${\cal C}$ is finitely generated and projective as
a left $A$-module.\\
Applying the results of the previous Section, we find a Morita
context connecting $B$ and ${}^*{\cal C}$. We will now show that
there is another Morita context connecting $B'$ and ${}^*{\cal C}$,
and that there is a morphism between the two Morita contexts.
We already know that $A$ is a $(B,{}^*{\cal C})$-bimodule, and this
implies that it is also a $(B',{}^*{\cal C})$-bimodule. We also have

\begin{lemma}\lelabel{2.2}
$Q'$ is a $({}^*{\cal C},B')$-bimodule.
\end{lemma}

\begin{proof}
Since we know that $Q$ is a $({}^*{\cal C},B)$-bimodule,
it suffices to show that
$Q$ is a left ideal in ${}^*{\cal C}$. 
For all $f\in {}^*{\cal C}$,
$q\in Q$ and $c\in {\cal C}$, we have
\begin{eqnarray*}
&&\hspace*{-2cm}
((f\# q)\# I_{\cal C})(c)=(f\#(q\# I_{\cal C}))(c)\\
&=& (f\# (\rho^l\circ q))(c)\\
&=& (\rho^l\circ q)(c_{(1)}f(c_{(2)}))\\
&=& q(c_{(1)}f(c_{(2)}))x\\
&=& (f\# q)(c)x
\end{eqnarray*}
\end{proof}

Now we define maps
\begin{equation}\eqlabel{2.3.1}
\mu':\ Q'\ot_{B'} A\to {}^*{\cal C}~~;~~\mu'(q\ot_{B'} a)=q\# i(a)
\end{equation}
\begin{equation}\eqlabel{2.3.2}
\tau':\ A\ot_{{}^*{\cal C}} Q'\to B'~~;~~\tau'(a\ot_{{}^*{\cal C}}q)=
a\cdot q=q(xa)
\end{equation}
It is clear that $\mu'$ is well-defined. $\tau'$ is also well-defined:
for all $f\in {}^*{\cal C}$, $a\in A$ and $q\in Q'$, we have
$$ \tau'(a\ot_{{}^*{\cal C}}(f\# q))=a\cdot (f\# q)=
(a\cdot f)\cdot q=\tau'(a\cdot f)\ot_{{}^*{\cal C}} q$$
and for all $q\in Q'$ and $a\in A$, we have that
\begin{equation}
q(xa)\in B'
\end{equation}
since
$$xq(xa)=q(xa)x$$

\begin{theorem}\thlabel{2.3}
With notation as above, $(B',{}^*{\cal C}, B', A, Q',\tau',\mu')$ is
a Morita context, and we have a morphism of Morita contexts
$$(B',{}^*{\cal C}, B', A, Q',\tau',\mu')\to
(B,{}^*{\cal C}, B, A, Q,\tau,\mu)$$
\end{theorem}

\begin{proof}
We have to show that the following two diagrams are commutative.
$$
\begin{diagram}
A\ot_{{}^*{\cal C}}Q'\ot_{B'} A&\rTo{\tau'\ot I_A}&B'\ot_{B'} A\\
\dTo^{I_A\ot \mu'}&&\dTo^{\cong}\\
A\ot_{{}^*{\cal C}}{}^*{\cal C}&\rTo^{\cong}&A
\end{diagram}
~~~~~~~~~~
\begin{diagram}
Q'\ot_{B'}A\ot_{{}^*{\cal C}}Q'&\rTo{\mu\ot I_{Q'}}&{}^*{\cal
C}\ot_{{}^*{\cal C}}Q'\\
\dTo^{I_{Q'}\ot \tau'}&&\dTo^{\cong}\\
Q'\ot_{B'}B'&\rTo^{\cong}&Q'
\end{diagram}
$$
Take $a,a'\in A$ and $q,q'\in Q'$. We compute
\begin{eqnarray*}
&&\hspace*{-2cm}
a\cdot\mu(q\ot a')=a\cdot(q\#i(a'))\\
&=& (q\# i(a'))(xa)=q(xa)a'= \tau(a\ot q)a'
\end{eqnarray*}
and this proves that the first diagram commutes. The second diagram
commutes if
$$(q\#i(a))\# q'=q\#i(q'(xa))$$
Indeed, for all $c\in {\cal C}$,
\begin{eqnarray*}
&&\hspace*{-2cm}
(q\#i(q'(xa)))(c)=q(c)q'(xa)=q'(q(c)xa)\\
&=& q'(c_{(1)}q(c_{(2)})a)=((q\#i(a))\# q')(c)
\end{eqnarray*}
The second statement is obvious: the morphism is given by
the inclusion maps $B'\subset B$, $Q'\subset Q$, and the identity
maps on $A$ and ${}^*{\cal C}$.
\end{proof}

We now present the coinvariants version of \prref{1a.2},
giving necessary and sufficient conditions for $\tau'$
to be surjective. It will follow
that our two Morita contexts coincide if $\tau'$ is surjective.

\begin{theorem}\thlabel{2.4}
Consider the Morita context $(B',{}^*{\cal C}, A, Q',\tau',\mu')$
of \thref{2.3}. The following statements are equivalent:\\
1) $\tau'$ is surjective (and, a fortiori, bijective);\\
2) there exists $\Lambda\in Q'$ such that $\Lambda(x)=1$;\\
3) for every right ${}^*{\cal C}$-module $M$, the map
$$\omega_M:\ M\ot_{{}^*{\cal C}} Q'\to M^{{}^*{\cal C}},~~
\omega_M(m\ot_{{}^*{\cal C}} q)=m\cdot q$$
is bijective.
\end{theorem}

\begin{proof}
$\ul{1)\Rightarrow 2)}$. If $\tau'$ is surjective, then there exist
$a_j\in A$ and $q_j\in Q'$ such that
$$1=\tau(\sum_j a_j\ot_{{}^*{\cal C}}q_j)=\sum_j q_j(xa_j)=
\sum_j(i(a_j)\# q_j)(x)$$
$\Lambda=\sum_i(a_j)\# q_j\in Q'$ because $Q'$ is a left ideal
in ${}^*{\cal C}$.\\
$\ul{2)\Rightarrow 3)}$. Define $\eta_M:\ M^{{}^*{\cal C}}\to 
M\ot_{{}^*{\cal C}} Q'$ as follows:
$$\eta_M(m)=m\ot_{{}^*{\cal C}}\Lambda$$
It is clear that $\omega_M\circ \eta_M=I_{M^{{}^*{\cal C}}}$. Furthermore,
for all $m\in M$ and $q\in Q'$,
$$\eta_M(\omega_M(m\ot_{{}^*{\cal C}}q))=
(m\cdot q)\ot_{{}^*{\cal C}}\Lambda=
m\ot_{{}^*{\cal C}}q\#\Lambda =m\ot_{{}^*{\cal C}}q$$
since
$$q\#\Lambda=q\Lambda(x)=q$$
$\ul{3)\Rightarrow 1)}$. $\omega_A=\tau'$ is bijective.
\end{proof}

\begin{theorem}\thlabel{2.5}
Consider the Morita context $(B',{}^*{\cal C},  A, Q',\tau',\mu')$
of \thref{2.3}. Assume that $\tau'$ is surjective, and take 
$\Lambda\in Q'$ such that $\Lambda(x)=1$. Then we have
the following properties.\\
1) $A$ and $Q$ are generators as left, resp. right, $B$-modules.\\
2) $A$ and $Q$ are finitely generated projective as right, resp. left,
${}^*{\cal C}$-modules.\\
3) We have bimodule isomorphisms
$$f:\ A\to {}_{{}^*{\cal C}}\Hom(Q,{}^*{\cal C});~~
f(a)(q)=q\# i(a)$$
$$g:\ Q\to \Hom_{{}^*{\cal C}}(A,{}^*{\cal C});~~
g(q)(a)=q\# i(a)$$
4) We have algebra isomorphisms
$$k:\ B\to {}_{{}^*{\cal C}}\End(Q);~~k(b)(q)=q\# i(b)$$
$$l: B\to \End(A)_{{}^*{\cal C}};~~l(b)(a)=ba$$
5) For all $M\in {\cal M}^{\cal C}$, $M^{{}^*{\cal C}}=M^{{\rm co}{\cal C}}$.
In particular $B=B'$.\\
6) $Q=Q'$.\\
7) The two Morita contexts in \thref{2.4} coincide.\\
8) $\Lambda\#\Lambda=\Lambda$ and $\Lambda\#{}^*{\cal C}\#\Lambda=\Lambda\#
B\cong B$.\\ 
9) For all $V\in {\cal M}_B$, the map
$$\eta_V:\ V\to (V\ot_B A)^{{\rm co}{\cal C}};~~\eta_V(v)=v\ot_B 1$$
is an isomorphism.\\
10) $B$ is a $B$-direct summand of $A$.\\
\end{theorem}

\begin{proof}
1), 2), 3) and 4) follow immediately from the Morita Theorems,
see \cite[II.3.4]{Ba}.\\
5) From \thref{2.4}, we know that there exists $\Lambda\in Q$
such that $\Lambda(x)=1$. Take $m\in M^{{}^*{\cal C}}$. Then
$$m=m\Lambda(x)=m\cdot\Lambda=m_{[0]}\Lambda(m_{[1]})$$
and
\begin{eqnarray*}
\rho(m)&=&\rho(m_{[0]}\Lambda(m_{[1]}))=m_{[0]}\ot_A m_{[1]}\Lambda(m_{[2]})\\
&=& m_{[0]}\ot_A\Lambda(m_{[1]})x=m_{[0]}\Lambda(m_{[1]})\ot_Ax=m\ot_Ax
\end{eqnarray*}
so it follows that $m\in M^{{\rm co}{\cal C}}$.\\
6) Look at the commutative diagram
$$\begin{diagram}
A\ot_{{}^*{\cal C}} Q'&\rTo^{\tau'}&B'\\
\dTo&&\dTo^{=}\\
A\ot_{{}^*{\cal C}} Q&\rTo^{\tau}&B
\end{diagram}$$
From the fact that $B=B'$ and $\tau'$ is surjective, we easily
deduce that $\tau$ is surjective. Applying 3) and its corresponding
property in \prref{1a.2a}, we find
$$Q\cong \Hom_{{}^*{\cal C}}(A,{}^*{\cal C})\cong Q'$$
7) now follows immediately from 5) and 6), and 8), 9) and 10)
follow from the corresponding properties in \prref{1a.2}.
\end{proof}

Now let us look at the map $\mu$. 
If we assume that
${\cal C}$ is finitely generated and projective as a left $A$-module.
As we already noticed, this implies that ${\cal M}^{\cal C}\cong
{\cal M}_{{}^*{\cal C}}$, the two Morita contexts coincide, and
we can apply \prref{1a.2}. Let us state the result, for completeness
sake. From \cite[II.3.4]{Ba}, recall that
we have ring morphisms
$$\pi:\ {}^*{\cal C}\to {}_B\End(A)^{\rm op},~~\pi(f)(a)=f(xa)$$
$$\pi':\ {}^*{\cal C}\to \End_B(Q),~~\pi'(f)(q)=f\#q$$
In fact $\pi={}^*{\rm can}$, cf. \prref{1.1}. We also have a
$({}^*{\cal C},B)$-bimodule map
$$\kappa:\ Q\to {}_B\Hom(A,B),~~\kappa(q)(a)=q(xa)$$
and a $(B,{}^*{\cal C})$-bimodule map
$$\kappa':\ A\to \Hom_B(Q,B),~~\kappa'(a)(q)=q(xa)$$
If $\mu$ is surjective, then $\pi$, $\pi'$, $\kappa$ and
$\kappa'$ are isomorphisms, and $A$ and $Q$ are finitely generated
and projective, respectively as left and a right $B$-module.
We now state some necessary and sufficient conditions for $\mu$
to be surjective.

\begin{theorem}\thlabel{2.6}
Assume that 
${\cal C}$ is finitely generated and projective as a left $A$-module,
and
consider the Morita context $(B=B',{}^*{\cal C}, A, Q=Q',\tau=\tau',\mu
=\mu')$ of \thref{2.3}. 
Then the following assertions are equivalent.\\
1) $\mu:\ Q\ot_B A\to {}^*{\cal C}$ is surjective (and, a fortiori,
bijective);\\
2) $({}^*{\cal C},x)$ satisfies the Weak Structure Theorem;
3) $A$ is a right ${}^*{\cal C}$-generator;\\
4) $A$ is projective as a left $B$-module and $\pi$ is bijective;\\
5) $A$ is projective as a left $B$-module, $\pi$ is injective,
and $\kappa$ is surjective;\\
6) $Q$ is projective as a right $B$-module, $\pi'$ is injective,
and $\kappa'$ is surjective;\\
7) $A$ is projective as a left $B$-module and $({\cal C},x)$ is a Galois
coring.
\end{theorem}

\begin{proof}
The equivalence of 1)-6) follows immediately from \prref{1a.2}.
$\ul{4)\Leftrightarrow 7)}$ follows from \prref{1.1}, using the fact that
 ${\cal C}$
is finitely generated and projective as left $A$-module.
\end{proof}

\section{Cleft entwining structures}\selabel{3}
In this Section, we look at the particular situation where
${\cal C}=A\ot C$ arises from an entwining structure $(A,C,\psi)$.
First observe that
$${}^*{\cal C}={}_A\Hom(A\ot C,A)\cong \Hom(C,A)$$
as a $k$-module. The ring structure on ${}^*{\cal C}$ induces a
$k$-algebra structure on $\Hom(C,A)$, and this $k$-algebra is
denoted $\#(C,A)$. The product is given by the formula
\begin{equation}\eqlabel{3.1}
(f\# g)(c)=f(c_{(2)})_{\psi}g(c_{(1)}^{\psi})
\end{equation}
We have a natural algebra homomorphism $i:\ A\to \#(C,A)$
given by
\begin{equation}\eqlabel{3.2}
i(a)(c)=\varepsilon_C(c)a
\end{equation}
and we have, for all $a\in A$ and $f:\ C\to A$:
\begin{equation}\eqlabel{3.3}
(i(a)\# f)(c)=a_{\psi}f(c^{\psi})~~{\rm and}~~
(f\#i(a))(c)=f(c)a
\end{equation}
$R=\Hom(C,A)$ will denote the $k$-algebra with the usual convolution
product, that is
\begin{equation}\eqlabel{3.4}
(f* g)(c)=f(c_{(1)})g(c_{(2)})
\end{equation}
The fact that we have two multiplications on $\Hom(C,A)$, namely the
usual convolution $*$ and the smash product $\#$ makes the difference
between the general coring theory and the theory of entwined modules.\\
We fix a grouplike element $x\in G(C)$. Then $1\ot x\in G({\cal C})$,
and the results of \seref{1} can be applied to this situation.
The following are then easily verified:\\
$A\in {\cal M}(\psi)_A^C$. The right $C$-coaction is
\begin{equation}\eqlabel{3.5}
\rho^r(a)=a_{\psi}\ot x^{\psi}
\end{equation}
The ring of coinvariants is
\begin{equation}\eqlabel{3.6}
B'=A^{{\rm co}C}=\{b\in A~|~b_{\psi}\ot x^{\psi}=b\ot x\}
\end{equation}
The bimodule $Q'$ is naturally isomorphic
to
$$Q'=\{q\in \#(C,A)~|~q(c_{(2)})_{\psi}\ot c_{(1)}^{\psi}=q(c)\ot x\}$$
We have maps
$$\mu':\ Q'\ot_{B'} A\to\#(C,A),~~\mu'(q\ot_{B'} a)(c)=q(c)a$$
$$\tau':\ A\ot_{\#(C,A)} Q'\to B',~~\tau'(a\ot q)=a_{\psi}q(x^{\psi})$$
and $(B,\#(C,A),A,Q,\tau',\mu')$ is a Morita context.\\
For $M\in {\cal M}_{\#(C,A)}$, the module of invariants is given by
$$M^R=\{m\in M~|~mf=mf(x),~{\rm for~all~}f\in R=\#(C,A)\}$$
From \thref{2.4}, we obtain immediately:

\begin{theorem}\thlabel{3.0a}
With notation as above, the following assertions are equivalent:\\
1) $\tau'$ is surjective;\\
2) there exists a $\Lambda\in Q$ such that $\Lambda(x)=1$;\\
3) for all $M\in {\cal M}_{\#(C,A)}$, the map
$$\omega_M:\ M\ot_{\#(C,A)} Q\to M^R,~~\omega_M(m\ot q)=m\cdot q$$
is bijective.
\end{theorem}

Now assume that $C$ is finitely generated projective as a $k$-module,
and let $\{c_j,c_j^*\}$ be a finite dual basis. Then we have a
natural isomorphism
$$\Hom(C,A)\cong A\ot C^*$$
The multiplication $\#$ on $\Hom(C,A)$ can be translated into a
multiplication on $A\ot C^*$. The $k$-algebra that we obtain in this
way is denoted $A\# C^*$. The multiplication can be described as 
follows (cf. e.g. \cite[Sec. 2.3]{CMZ}): define $R:\ C^*\ot A\to A\ot C^*$
by
$$R(c^*\ot a)=a_R\ot c^*_R=\sum_j \lan c^*,c_j^{\psi}\ran a_{\psi}\ot c_j^*$$
Then
$$(a\# c^*)(b\# d^*)=a_Rb\# (d_R^**c^*)$$
We have maps
$$\mu:\ Q\ot_B A\to A\# C^*,~~\mu(q\ot_B a)=\sum_j\lan q,c_j\ran a\# c_j^*$$
$$\pi:\ A\# C^*\to \End_B(A),~~\pi(b\# c^*)(a)=\lan c^*,xa\ran b$$
$${\rm can}:\ A\ot_B A\to A\ot C,~~{\rm can}(a\ot b)=ab_{\psi}\ot c^{\psi}$$
For every $M\in {\cal M}(\psi)_A^C$, we have
$$\varepsilon_M:\ M^{{\rm co}C}\ot_B A\to M,~~\phi_M(m\ot a)=ma$$
From \thref{2.6}, we immediately obtain the following:

\begin{theorem}\thlabel{3.0b}
Let $(A,C,\psi)$ be an entwining structure, and $x\in C$ grouplike,
and assume that $C$ is finitely generated and projective as a $k$-module.
Then the following assertions are equivalent:\\
1) $\mu$ is surjective (and a fortiori bijective);\\
2) $\varepsilon_M$ is bijective, for every $M\in {\cal M}(\psi)_A^C$;\\
3) $A$ is a right $A\#C^*$-generator;\\
4) $A$ is projective as a left $A$-module, and $\pi$ is bijective;\\
5) $A$ is projective as a left $A$-module, and ${\rm can}$ is bijective,
i.e. $A$ is a $C$-coalgebra Galois extension in the sense of
\cite{Br1}.
\end{theorem}

\begin{proposition}\prlabel{3.1}
Assume that $\lambda:\ C\to A$ is convolution invertible, with
convolution inverse $\lambda^{-1}$. Then the following assertions are
equivalent:\\
1) $\lambda\in Q'$;\\
2) for all $c\in C$, we have
\begin{equation}\eqlabel{3.1.1}
\lambda^{-1}(c_{(1)})\lambda(c_{(3)})_{\psi}\ot c_{(2)}^{\psi}=
\varepsilon(c)1_A\ot x
\end{equation}
3) for all $c\in C$, we have
\begin{equation}\eqlabel{3.1.2}
\lambda^{-1}(c_{(1)})\ot c_{(2)}=\lambda^{-1}(c)_{\psi}\ot x^{\psi}
\end{equation}
Notice that condition 3) means that $\lambda^{-1}$ is right
$C$-colinear. If such a $\lambda\in Q'$ exists, then we call
$(A,C,\psi,x)$ cleft.
\end{proposition}

\begin{proof}
$\ul{1) \Rightarrow 2)}$.
$$\lambda^{-1}(c_{(1)})\lambda(c_{(3)})_{\psi}\ot c_{(2)}^{\psi}=
\lambda^{-1}(c_{(1)})\lambda(c_{(2)})\ot x=
\varepsilon_C(c)1_A\ot x$$
$\ul{2) \Rightarrow 3)}$.
\begin{eqnarray*}
&&\hspace*{-2cm}
\lambda^{-1}(c)_{\psi}\ot x^{\psi}=
\varepsilon(c_{(1)})1_A\lambda^{-1}(c_{(2)})_{\psi}\ot x^{\psi}\\
{\rm \eqref{3.1.1}}~~~~&=& \lambda^{-1}(c_{(1)})\lambda(c_{(3)})_{\Psi}
\lambda^{-1}(c_{(4)})_{\psi}\ot c_{(2)}^{\Psi\psi}\\
{\rm \eqref{2.1.1.1}}~~~~&=&\lambda^{-1}(c_{(1)})
\Bigl(\lambda(c_{(3)})\lambda^{-1}(c_{(4)})\Bigr)_{\psi}\ot c_{(2)}^{\psi}\\
&=& \lambda^{-1}(c_{(1)})\ot c_{(2)}
\end{eqnarray*}
$\ul{3) \Rightarrow 1)}$.
\begin{eqnarray*}
&&\hspace*{-2cm} \lambda(c_{(2)})_{\psi}\ot c_{(1)}^{\psi}=
\lambda(c_{(1)})\lambda^{-1}(c_{(2)})
\lambda(c_{(4)})_{\psi}\ot c_{(3)}^{\psi}\\
{\rm \eqref{3.1.2}}~~~~&=& 
\lambda(c_{(1)})\lambda^{-1}(c_{(2)})_{\Psi}
\lambda(c_{(3)})_{\psi}\ot x^{\Psi\psi}\\
{\rm \eqref{2.1.1.1}}~~~~&=&
\lambda(c_{(1)})\Bigl(\lambda^{-1}(c_{(2)})
\lambda(c_{(3)})\Bigr)_{\psi}\ot x^{\psi}\\
&=& \lambda(c)\ot x
\end{eqnarray*}
\end{proof}

\begin{proposition}\prlabel{3.2}
Assume that $(A,C,\psi,x)$ is a cleft entwining structure. Then the
map $\tau'$ in the associated Morita context is surjective.
\end{proposition}

\begin{proof}
Let $\lambda$ be as in \prref{3.1}. From condition 3) in \prref{3.1},
we deduce that
$$\lambda^{-1}(x)\ot x=\lambda^{-1}(x)_{\psi}\ot x^{\psi}$$
hence $\lambda^{-1}(x)\in B'$, and $\Lambda=\lambda\# i(\lambda^{-1}(x))
\in Q'$, since $Q'$ is a right $B'$-module. Now
$$\Lambda(x)=\lambda(x)\lambda^{-1}(x)=\varepsilon_C(x)=1$$
and it follows from \thref{3.0a} that $\tau$ is surjective.
\end{proof}

We say that the entwining structure $(A,C,\psi,x)$ satisfies the {\sl right
normal basis property} if there exists a left $B'$-linear and
right $C$-colinear isomorphism
$B'\ot C\to A$. $(A,C,\psi,x)$ satisfies the Strong (resp. Weak)
Structure Theorem if $(A\ot C,1\ot x)$ 
satisfies the Strong (resp. Weak) Structure Theorem.
We can now state
our main result. 

\begin{theorem}\thlabel{3.3}
Let $(A,C,\psi,x)$ be an entwining structure with a fixed grouplike element.
The following assertions are equivalent:\\
1) $(A,C,\psi,x)$ is cleft;\\
2) $(A,C,\psi,x)$ satisfies the
Strong Structure Theorem and the right normal basis property ;\\
3) $(A,C,\psi,x)$ is Galois, and satisfies the right normal basis property;\\
4) the map ${}^*{\rm can}:\ \#(C,A)\to \End_B(A)^{\rm op}$ is
bijective and $(A,C,\psi,x)$ satisfies the
right normal basis property.
\end{theorem}

\begin{proof}
$\ul{1)\Rightarrow 2)}$. We take $\lambda\in Q'$ as in \prref{3.1},
and $M\in {\cal M}(\psi)_A^C$. For any $m\in M$, we have
\begin{eqnarray*}
\rho(m\cdot \lambda)&=&
\rho(m_{[0]}\lambda(m_{[1]}))=m_{[0]}\lambda(m_{[2]})_{\psi}
\ot \lambda(m_{[1]})^{\psi}\\
&=& m_{[0]}\lambda(m_{[1]})\ot x=m\cdot \lambda \ot x
\end{eqnarray*}
hence $m\cdot \lambda\in M^{{\rm co}C}$, and we have a
well-defined map
$$\gamma_M:\ M\to M^{{\rm co}C}\ot_{B'} A,~~\gamma_M(m)=
m_{[0]}\cdot\lambda\ot_{B'}\lambda^{-1}(m_{[1]})$$
and we compute easily that
$$\varepsilon_M(\gamma_M(m))=m_{[0]}\cdot (\lambda\# i(\lambda^{-1}(m_{[1]})))
=m_{[0]}\lambda(m_{[1]})\lambda^{-1}(m_{[2]})=m$$
Recall that
\begin{equation}\eqlabel{3.3.1}
\tau(a\ot\lambda)=a_{\psi}\lambda(x^{\psi})\in B'
\end{equation}
Take
$a\in A$ and $m\in M^{{\rm co}C}$. Then
\begin{eqnarray*}
&&\hspace*{-2cm}
\gamma_M(\varepsilon_M(m\ot_{B'} a))=\gamma_M(ma)=
pa_{\psi}\cdot \lambda\ot_{B'} \Lambda^{-1}(x^{\psi})\\
&=& pa_{\psi\Psi}\lambda(x^{\Psi})\ot_{B'} \Lambda^{-1}(x^{\psi})\\
{\rm \eqref{3.3.1}}~~~~
&=& p\ot_{B'}a_{\psi\Psi}\lambda(x^{\Psi}) \Lambda^{-1}(x^{\psi})\\
{\rm \eqref{2.1.1.3}}~~~~
&=& p\ot_{B'} a_{\psi}\lambda((x^{\psi})_{(1)})\lambda^{-1}((x^{\psi})_{(2)})\\
&=& p\ot_{B'} a_{\psi}\varepsilon(x^{\psi})\\
{\rm \eqref{2.1.1.4}}~~~~
&=& p\ot_{B'} a
\end{eqnarray*}
and this proves that $(A,C,\psi,x)$ satisfies the Weak Structure
Theorem. From \prref{3.2}, we know that the map $\tau'$ in the
Morita context is surjective, and part 3) of \thref{2.5} tells us
that $(A,C,\psi,x)$ satisfies the Strong Structure Theorem.\\
Take $M\in {\cal M}(\psi)_A^C$, and consider the maps
$$k:\ M\to M^{{\rm co}C}\ot C,~~k(m)=m_{[0]}\cdot\lambda\ot m_{[1]}=
m_{[0]}\lambda(m_{[1]})\ot m_{[2]}$$
$$k^{-1}:\ M^{{\rm co}C}\ot C\to M,~~k^{-1}(m\ot c)=m\lambda^{-1}(c)$$
It is clear that $k^{-1}(k(m))=m$, for all $m\in M$.
For $m\in M^{{\rm co}C}$. Then $\rho(ma)=ma_{\psi}\ot x^{\psi}$,
and we compute
\begin{eqnarray*}
&&\hspace*{-2cm}
k(k^{-1}(m\ot c))=(m\lambda^{-1}(c)_{\psi}\cdot\lambda\ot x^{\psi}\\
{\rm \eqref{3.1.2}}~~~~
&=& (m\lambda^{-1}(c_{(1)})\cdot\lambda\ot c_{(2)}\\
&=& m\lambda^{-1}(c_{(1)})_{\psi}\lambda(x^{\psi})\ot c_{(2)}\\
{\rm \eqref{3.1.2}}~~~~
&=& m\lambda^{-1}(c_{(1)}) \lambda(c_{(2)})\ot c_{(3)}\\
&=& m\ot c
\end{eqnarray*}
It is obvious that $k$ is right $C$-colinear. Now $A\in {\cal M}(\psi)_A^C$,
so we find a right $C$-colinear isomorphism $A\cong B\ot C$. It is
also left $B$-linear, since the right $C$-coaction on $A$ is left
$B$-linear.\\
$\ul{2)\Rightarrow 3)}$ follows from part 4) of \prref{1.1}.\\
$\ul{3)\Rightarrow 4)}$ follows from part 1) of \prref{1.1}.\\
$\ul{4)\Rightarrow 1)}$.
From the right normal basis property, we know that there exists
a left $B'$-linear, right $C$-colinear isomorphism
$h:\ B'\ot C\to A$. We consider the maps
$$\lambda:\ C\to A,~~\lambda(c)=h(1\ot c)$$
and
$$j=(I_{B'}\ot \varepsilon_C)\circ h^{-1}:\ A\to B'$$
It is then clear that $h$ is right $C$-colinear and $j$ is
left $B'$-linear. Take $a\in A$, and write
$$h^{-1}(a)=\sum_i b_i\ot c_i$$
Then
\begin{eqnarray*}
&&\hspace*{-2cm}\sum_i b_ih(1\ot c_{i(1)})\ot c_{i(2)}=
\sum_i h(b_i\ot c_{i(1)})\ot c_{i(2)}\\
&=& (h\ot I_C)\rho(\sum_i b_i\ot c_i)\\
&=& \rho(h(\sum_i b_i\ot c_i))=\rho(a)=a_{\psi}\ot x^{\psi}
\end{eqnarray*}
Apply $j\ot I_C$ to both sides:
\begin{eqnarray*}
&&\hspace*{-2cm} j(a_{\psi})\ot x^{\psi}=
\sum_i b_i(j\circ h)(1\ot c_{i(1)})\ot c_{i(1)}\\
&=& \sum_i b_i(I_B\ot \varepsilon_C)(1\ot c_{i(1)})\ot c_{i(1)}\\
&=& \sum_i b_i\ot c_i=h^{-1}(a)
\end{eqnarray*}
Now let $q=({}^*{\rm can})^{-1}(j)$. We are done if we can show
that $\lambda$ is a convolution inverse of $q$, by \prref{3.1}.
The fact that $\lambda$ is right $C$-colinear means
\begin{equation}\eqlabel{3.3.2}
\lambda(c_{(1)})\ot c_{(2)}=\lambda(c)_{\psi}\ot x^{\psi}
\end{equation}
and we compute, for all $c\in C$,
\begin{eqnarray*}
&&\hspace*{-2cm}
(\lambda * q)(c)=\lambda(c_{(1)})q(c_{(2)})=
\lambda(c)_{\psi}q( x^{\psi})\\
&=&{}^*{\rm can}(q)(\lambda(c))=j(\lambda(c))\\
&=& ((I_{B'}\ot \varepsilon_C)\circ h^{-1}\circ h)(1\ot c)=\varepsilon_C(c)1_A
\end{eqnarray*}
as needed. For all $a\in A$, we have
\begin{eqnarray*}
&&\hspace*{-2cm}
{}^*{\rm can}(q*\lambda)(a)=a_{\psi}(q*\lambda)(x^{\psi})\\
&=& a_{\psi}q((x^{\psi})_{(1)})\lambda((x^{\psi})_{(2)})\\
{\rm \eqref{2.1.1.3}}~~~~
&=& a_{\psi\Psi}q(x^{\Psi})\lambda(x^{\psi})\\
&=& {}^*{\rm can}(q)(a_{\psi})\lambda(x^{\psi})\\
&=& j(a_{\psi})\lambda(x^{\psi})\\
&=& j(a_{\psi})h(1\ot x^{\psi})\\
&=& h(j(a_{\psi})\ot x^{\psi})=h(h^{-1}(a))=a
\end{eqnarray*}
This proves that
$${}^*{\rm can}(q*\lambda)=I_A={}^*{\rm can}(\eta_A\circ\varepsilon_C)$$
and
$$q*\lambda=\eta_A\circ\varepsilon_C$$
by the injectivity of ${}^*{\rm can}$. This finishes the proof of
the Theorem.
\end{proof}

\section{Factorization structures and the CFM Morita context}
\selabel{5}
Let $(A,S,R)$ be a factorization structure, and consider the
smash product $R=A\#_R S$. We fix an algebra map $\chi:\ S\to k$.
Then the map
$$X:\ R=A\#_R S\to A,~~X(a\# s)=\chi(s)a$$
satisfies the conditions of \seref{1a} (with right replaced by left):
$X$ is left $A$-linear, $X(rX(s))=X(rs)$, and $X(1)=1$.
We can therefore apply the results of \seref{1a}. In particular,
we obtain that $A$ is a left $R$-module:
$$(a\# s)\rightact b=X((a\# s)b)=X(ab_R\# s_R)=\chi(s_R)ab_R$$
and $b\in B=A^R$ if and only if
$\chi(s_R)b_R=\chi(s)b$ for all $s\in S$. Also
$\sum_i a_i\# s_i\in Q$ if and only if
$$\sum_i a_{iR}\#t_Rs_i=\chi(t)\sum_i a_{i}\#s_i$$
for all $t\in S$. We have a Morita context
$(B,A\#S,A,Q,\tau,\mu)$ with
$$\mu:\ A\ot_B Q\to A\# S,~~\tau(a\ot_B(\sum_i a_i\# s_i))=
\sum_i aa_i\# s_i$$
$$\tau:\ Q\ot_R A\to B,~~\tau(\sum_i a_i\# s_i)\ot_R a)=
\sum_i a_ia_R\chi(s_{iR})$$
Now we consider the following particular situation: $S=H$ is a 
bialgebra, $A$ is a left $H$-module algebra, and
$$R:\ H\ot A\to A\ot H,~~R(h\ot a)=h_{(1)}\cdot b\ot h_{(2)}$$
We also take $\chi=\varepsilon_H$. The above formulas take the
following form:
$$(a\#h)\rightact b=a(h\cdot b)$$
$$b\in B~~\Longleftrightarrow~~h\cdot b=\varepsilon(h)b$$
$$\sum_i a_i\# h_i\in Q~~\Longleftrightarrow~~
\sum_i h_{(1)}\cdot a_i\# h_{(2)}h_i=\varepsilon(h)\sum_i a_i\# h_i$$
for all $h\in H$.\\
In the particular situation where $H$ is a finite dimensional Hopf
algebra over a field $k$, there exists another Morita context
connecting $B$ and $A\# H$, due to Cohen, Fischman and Montgomery
(see \cite{CFM}).
The construction can be generalized to the case where
$H$ is a Frobenius Hopf algebra
over a commutative ring $k$ (see \cite{CF}).
This Morita context can be described as follows. Take a free
generator $t$ of the space of left integrals in $H$, and let $\lambda$
be the distinguished grouplike element in $H^*$:
$$ht=\varepsilon(h)t,~{\rm and}~th=\lambda (h)t$$
for all $h\in H$. Then $\lambda$ is an algebra map, and $A$ is
a $(B,A\# H)$-bimodule, the right $A\# H$-action is given by
$$a\rightact (b\# h)=\lambda(h_{(2)})\ol{S}(h_{(1)})\cdot(ab)$$
and we have a Morita context
$$(B,A\# H,A,A,\tilde{\tau},\tilde{\mu})$$
with
$$\tilde{\tau}:\ A\ot_RA\to B,~~\tilde{\tau}(a\ot_R b)=t\cdot (ab)$$
$$\tilde{\mu}:\ A\ot_BA\to R,~~\tilde{\mu}(a\ot_B b)=
a(t_{(1)}\cdot b)\# t_{(2)}$$
We refer to \cite{CFM} for the details. We will now show that
this Morita context can be obtained using \prref{1a.1} and
\thref{1a.4}.\\
If $H$ is Frobenius, then there exists a left integral $\varphi$
in $H^*$ such that $\lan \varphi,t\ran=1$. $\varphi$ is a free generator
of the space of left integrals in $H^*$, and
$(t_{(2)}\ot \ol{S}(t_{(1)}),\varphi)$ is a Frobenius system for
$H/k$ (see for example \cite[Theorem 31]{CMZ}). This means that
\begin{equation}\eqlabel{5.1.1}
ht_{(2)}\ot\ol{S}(t_{(1)})=t_{(2)}\ot\ol{S}(t_{(1)})h~~{\rm and}~~
\lan \varphi, t_{(2)}\ran\ol{S}(t_{(1)})=t_{(2)}\lan
\varphi,\ol{S}(t_{(1)})\ran =1
\end{equation}
for all $h\in H$.

\begin{proposition}\prlabel{5.1}
Let $H$ be a Frobenius Hopf algebra, let $t$ and $\varphi$ be
as above, and take a left $H$-module algebra $A$. 
Then $A\# H/A$ is Frobenius, with Frobenius system
$$\Bigl(e=(1\# t_{(2)})\ot_A(1\# \ol{S}(t_{(1)})),
\ol{\nu}=I_A\#\varphi\Bigr)$$
\end{proposition}

\begin{proof}
We first show that $e$ is a Casimir element. Indeed,
for all $a\in A$ and $h\in H$, we have
\begin{eqnarray*}
&&\hspace*{-2cm}
\Bigl((1\# t_{(2)})\ot_A(1\# \ol{S}(t_{(1)})\Bigr)(a\# h)\\
&=&(1\# t_{(3)})\ot_A \bigl(\ol{S}(t_{(2)})\cdot a\#
\ol{S}(t_{(1)})h\bigr)\\
&=& \bigl((t_{(3)}\ol{S}(t_{(2)}))\cdot a\# t_{(4)}\bigr)
\ot_A \bigl(1\# \ol{S}(t_{(1)})h\bigr)\\
&=& (a\# t_{(2)})\ot_A \bigl(1\# \ol{S}(t_{(1)})h\bigr)\\
&=& (a\# h)\Bigl((1\# t_{(2)})\ot_A(1\# \ol{S}(t_{(1)})\Bigr)
\end{eqnarray*}
It is obvious that $\ol{\nu}$ is left $A$-linear. It is also
right $A$-linear since
\begin{eqnarray*}
&&\hspace*{-2cm}
\ol{\nu}((1\# h)a)=\ol{\nu}(h_{(1)}a\# h_{(2)})\\
&=& \lan \varphi,h_{(2)}\ran h_{(1)}a=
\lan \varphi,h\ran a
\end{eqnarray*}
Finally, using \eqref{5.1.1}, we find that
\begin{eqnarray*}
&&\hspace*{-2cm}
\ol{\nu}(1\# t_{(2)})(1\# \ol{S}(t_{(1)}))=
1\#\ol{S}(\lan\varphi,t_{(2)}\ran t_{(1)})\\
&=&1\#\ol{S}(\lan\varphi,t\ran 1)=1\#1
\end{eqnarray*}
and
$$(1\# t_{(2)})\ol{\nu}(1\# \ol{S}(t_{(1)}))=
1\# t_{(2)}\lan\varphi, \ol{S}(t_{(1)})\ran=1\# 1$$
\end{proof}

\begin{corollary}\colabel{5.2}
As in \prref{5.1}, let $H$ be a Frobenius Hopf algebra,
and $A$ a left $H$-module algebra. Then $A$ and $Q$ are
isomorphic as $(A,A\# H)$-bimodules and the Morita contexts
from \prref{1a.1} and \cite{CFM},\cite{CF} are isomorphic.
\end{corollary}

\begin{proof}
The fact that $A$ and $Q$ are isomorphic follows immediately
from \thref{1a.4} and \prref{5.1}. The connecting isomorphisms
are
$$\alpha:\ A\to Q,~~\alpha(a)=t_{(1)}\cdot a\# t_{(2)}$$
and $\alpha=I_A\# \varphi_{|Q}$. Let us check that the right
$A\#H$-action on $A$ transported from the one on $Q$ coincides
with the $A\#H$-action from \cite{CFM}:
\begin{eqnarray*}
&&\hspace*{-2cm}
a\leftact (b\#h)=\ol{\nu}(\alpha(a)(b\# h))\\
&=& \ol{\nu}\Bigl(\bigl( (t_{(1)}\cdot a)\# t_{(2)}\bigr)
(b\# h)\Bigr)\\
&=& \ol{\nu}\Bigl((t_{(1)}\cdot a)(t_{(2)}\cdot b)\# t_{(3)}h\Bigr)\\
&=& \lan\varphi,t_{(2)}h\ran t_{(1)}\cdot (ab)\\
&=& \lan\varphi,t_{(2)}h_{(3)}\ran (t_{(1)}h_{(2)}\ol{S}(h_{(1)}))\cdot (ab)\\
&=& \lan\varphi,th_{(2)}\ran\ol{S}(h_{(1)})\cdot (ab)\\
&=& \lan\lambda h_{(2)}\ran \lan\varphi, t\ran \ol{S}(h_{(1)})\cdot (ab)\\
&=& \lan\lambda h_{(2)}\ran \ol{S}(h_{(1)})\cdot (ab)
\end{eqnarray*}
as needed.
\end{proof}

\section{Cleft factorization structures}\selabel{6}
As in the beginning of \seref{5}, let $(A,S,R)$ be a factorization
structure, and $\chi:\ S\to k$ an algebra map. Recall that
$q=\sum_i a_i\# s_i\in Q$ if and only if
\begin{equation}\eqlabel{6.1.1}
\sum_i a_{iR}\# t_Rs_i=\chi(t)\sum_i a_i\# s_i
\end{equation}
for all $t\in S$. Take $q=\sum_i a_i\# s_i\in Q$, and assume that
$q$ is invertible in $A^{\rm op}\ot S$, i.e. there exists
$\ol{q}=\sum_j \ol{a}_j\# \ol{s}_j\in A\# S$ such that
\begin{equation}\eqlabel{6.1.*}
\sum_{i,j}  a_i\ol{a}_j\# \ol{s}_js_i=
\sum_{i,j}  \ol{a}_ja_i\# s_i\ol{s}_j=1_A\#1_S
\end{equation}

\begin{proposition}\prlabel{6.1}
Let $q=\sum_i a_i\# s_i\in A\# S$ be invertible in
$A^{\rm op}\ot S$, with inverse $\ol{q}=\sum_j \ol{a}_j\# \ol{s}_j$.
Then the following assertions are equivalent:\\
1) $q\in Q$;\\
2) for all $t\in S$:
\begin{equation}\eqlabel{6.1.2}
\sum_{i,j} (a_j)_R\ol{a}_i\#\ol{s}_it_Rs_j=
\chi(t)1_A\#1_S
\end{equation}
3) for all $t\in S$:
\begin{equation}\eqlabel{6.1.3}
\sum_j \chi(t_R)(\ol{a}_j)_R\# \ol{s}_j=
\sum_j\ol{a}_j\#\ol{s}_jt
\end{equation}
In this situation, we call $(A,S,R,\chi)$ a cleft factorization
structure.
\end{proposition}

\begin{proof}
$\ul{1)\Rightarrow 2)}$. Using \eqref{6.1.1} and \eqref{6.1.*},
we find, for all $t\in S$:
$$\sum_{i,j} (a_j)_R\ol{a}_i\#\ol{s}_it_Rs_j=
\chi(t)\sum_{i,j}  a_j\ol{a}_i\# \ol{s}_is_j=
\chi(t)1_A\#1_S$$
$\ul{2)\Rightarrow 3)}$. For all $t\in S$, we compute
\begin{eqnarray*}
&&\hspace*{-2cm}\sum_j \chi(t_R)(\ol{a}_j)_R\# \ol{s}_j=
\sum_j\chi(t_R)(\ol{a}_j)_R1_A\# 1_S\ol{s}_j\\
{\rm \eqref{6.1.2}}~~~~
&=& \sum_{i,j,k}(\ol{a}_j)_R(a_k)_r\ol{a}_i\#
\ol{s}_it_{Rr}s_k\ol{s}_j\\
{\rm \eqref{2.1.2.4}}~~~~
&=& \sum_{i,j,k}(\ol{a}_ja_k)_R\ol{a}_i\#
\ol{s}_it_{R}s_k\ol{s}_j\\
{\rm \eqref{6.1.*}}~~~~
&=& \sum_i (1_A)_R\ol{a}_i\#\ol{s}_it_R1_S\\
{\rm \eqref{2.1.2.1}}~~~~
&=& \sum_i\ol{a}_i\#\ol{s}_it
\end{eqnarray*}
$\ul{3)\Rightarrow 1)}$.
\begin{eqnarray*}
&&\hspace*{-2cm}
\sum_i (a_i)_R\# t_Rs_i=
\sum_{i,j,k} (a_i)_R\ol{a}_ja_k\# s_k\ol{s}_jt_Rs_i\\
{\rm \eqref{6.1.3}}~~~~
&=& \sum_{i,j,k}\chi(t_{Rr})(a_i)_R(\ol{a}_j)_ra_k\#
s_k\ol{s}_js_i\\
{\rm \eqref{2.1.2.4}}~~~~
&=& \sum_{i,j,k}\chi(t_{R})(a_i\ol{a}_j)_Ra_k\#
s_k\ol{s}_js_i\\
{\rm \eqref{6.1.*}}~~~~
&=& \sum_k a_k\#s_k
\end{eqnarray*}
\end{proof}

\begin{proposition}\prlabel{6.2}
Assume that $(A,S,R,\chi)$ is cleft. Then we have an equivalence
of categories
$$F:\ {}_B{\cal M}\to{}_R{\cal M},~~F(N)=A\ot_BN$$
$$G:\ {}_R{\cal M}\to{}_B{\cal M},~~G(M)={}^RM$$
Consequently the map
$${\rm can}:\ A\ot_BA\to \Hom(S,A),~~{\rm can}(a\ot a')(s)=
a'_R\chi(s_R)a$$
is bijective.
\end{proposition}

\begin{proof}
1) We first prove that the functor $F$ is fully faithful.
This follows from \prref{1a.2a} after we show that the map
$\tau$ from the Morita context from \prref{1a.1} is surjective.
It suffices to show that there exists $q\in Q$ with
$(I_A\ot\chi)(q)=1$ (\prref{1a.2}).\\
Take $q\in Q$ as in \prref{6.1}. From \eqref{6.1.3} it follows
that
$$\sum_j\chi(t_R)\chi(\ol{s}_j)(\ol{a}_j)_R=
\sum_j\chi(t)\chi(\ol{s}_j)\ol{a}_j$$
and this means that $\sum_j\chi(\ol{s}_j)\ol{a}_j\in B$.
$Q$ is a left $B$-module, so
$$\sum_{i,j}\chi(\ol{s}_j)\ol{a}_ja_i\#s_i\in Q$$
and it follows from \eqref{6.1.*} that
$$(I_A\ot\chi)(\sum_{i,j}\chi(\ol{s}_j)\ol{a}_ja_i\#s_i)=
\sum_{i,j}\chi(\ol{s}_j)\chi(s_i)\ol{a}_ja_i=1$$
2) Now we show that $G$ is fully faithful, or, equivalently,
the counit of the adjunction $(F,G)$ is an isomorphism.
Recall that, for $M\in {}_R{\cal M}$:
$$\varepsilon_M:\ A\ot_B{}^RM\to M,~~
\varepsilon_M(a\ot m)=am$$
Take $q\in Q$ as in \prref{6.1}, and $m\in M$. Then
$q\cdot m\in {}^RM$ since
\begin{eqnarray*}
(1\# t)q\cdot m&=&\sum_i (1\# t)(q_i\# s_i)\cdot m\\
&=& (\sum_i (q_i)_R\#t_Rs_i)\cdot m\\
{\rm \eqref{6.1.1}}~~~~&=& \chi(t)q\cdot m
\end{eqnarray*}
Now define $\gamma_M:\ M\to A\ot_B{}^RM$ as follows:
$$\gamma_M(m)=\sum_j \ol{a}_j\ot_B q\ol{s}_jm$$
For all $m\in M$, we have, by \eqref{6.1.*}
$$\varepsilon_M(\gamma_M(m))=\sum_j \ol{a}_jq\ol{s}_jm=m$$
and, finally, for all $b\in A$ and $m\in {}^RM$:
\begin{eqnarray*}
&&\hspace*{-2cm}
\gamma_M(\varepsilon_M(b\ot_B m))=
\gamma_M(bm)=\sum_j \ol{a}_j\ot_B q\ol{s}_jbm\\
&=&\sum_j \ol{a}_j\ot_B X(q\ol{s}_jb)m=
\sum_j \ol{a}_j X(q\ol{s}_jb)\ot_Bm\\
&=& \sum_{i,j} \ol{a}_ja_i\chi(s_i\ol{s}_j)b\ot_B m=b\ot_B m
\end{eqnarray*}
We used the fact that $X(q\ol{s}_jb)\in {\rm Im}(\tau)\subset B$.
\end{proof}

\section*{Acknowledgement}
The authors thank D. Quinn and \c S. Raianu for helpful comments.

\end{document}